\documentclass[letterpaper, 10 pt, conference]{ieeeconf}  

\IEEEoverridecommandlockouts                              
\usepackage{graphics} 
\usepackage{epsfig} 
\usepackage{amsmath} 
\usepackage{amssymb}  

\usepackage[tight,footnotesize]{subfigure}

\usepackage{mathrsfs} 
\usepackage{mathrsfs} 
\usepackage{mathtools} 
\usepackage{hyperref} 
\usepackage[all]{hypcap} 
\usepackage{pdfsync}
\usepackage{bbm}
\usepackage{verbatim}

\hypersetup{
  colorlinks = true,
  urlcolor = black
}

\usepackage{graphics}
\graphicspath{{./gfx/}}
\DeclareGraphicsExtensions{.eps}
\usepackage{epstopdf}

\newcommand{\ignore}[1]{}

\newcommand{\uK}{\underline{K}}
\newcommand{\oK}{\overline{K}}

\newcommand{\oD}{\bar{\Delta}}
\newcommand{\uD}{\underline{\Delta}}

\newcommand{\Di}{\Delta^k}

\newcommand{\cA}{\mathcal{A}}

\newcommand{\cB}{\mathcal{B}}

\newcommand{\cP}{\mathcal{P}}

\newcommand{\ccab}{\mathcal{C}_{\alpha, \beta}}
\newcommand{\D}{\Delta}

\newcommand{\Fc}{{ \mathscr{F}}}

\newcommand{\cFt}{{ \mathscr{F}}_{t}}
\newcommand{\Gc}{{ \mathscr{G}}}

\newcommand{\calo}{{\mathcal{O}}}
\newcommand{\Exp}{{\sf E}}

\newcommand{\Pro}{{\sf P}}

\newcommand{\Hyp}{{\sf H}}

\newtheorem{theorem}{Theorem}[section]

\title{\LARGE \bf
Unstructured sequential testing in sensor networks
}

\author{Georgios Fellouris and Alexander Tartakovsky
\thanks{G. Fellouris is with the Department of Statistics, University of Illinois, Urbana-Champaign, 
119 Illini Hall, 725 South Wright Street, Champaign, IL 61820 USA {\tt \small fellouri@illinois.edu}.
A. Tartakovsky is with the  Department of Statistics, University of Connecticut, 
215 Glenbrook Road U-4120, Sorrs, CT 06269-4120, {\tt \small a.tartakov@uconn.edu}. This paper was originally submitted when both authors were with the Department of Mathematics at the University of Southern California.}}

\begin{document}

\thispagestyle{empty}
\pagestyle{empty}

\maketitle

\begin{abstract}
We consider the problem of quickly detecting a signal in a sensor network when the subset of sensors in which signal may be present is completely unknown. 
We formulate this problem as a sequential hypothesis testing problem with a simple null (signal is absent everywhere) and a composite alternative (signal is present somewhere).  We introduce a novel class of scalable sequential tests which, 
\textit{for any subset of affected sensors},  minimize the expected sample size for a decision asymptotically, that is as the error probabilities go to 0. 
Moreover, we propose sequential tests  that require minimal transmission activity from the sensors to the fusion center, 
while preserving this asymptotic optimality property. 
\end{abstract}

\section{Introduction}
\subsection{Problem formulation}
Consider $K$ sources of observations (sensors) which transmit their data to a global decision maker (fusion center). 
We assume that observations from different sensors are  independent and that, for each $1 \leq k \leq K$, sensor $k$ observes a sequence $(X^{k}_{t})_{t \in \mathbb{N}}$ of independent and identically distributed (i.i.d.) random variables with common density $f^{k}$ with respect to a  dominating, $\sigma$-finite measure $\nu^{k}$.  We denote by $\{\Fc_{t}^{k}, t\in \mathbb{N}\}$ the filtration generated by the observations at sensor $k$, i.e., 
$\Fc_{t}^{k}:=\sigma(X_{s}^{k}; 1 \leq s \leq t)$ for every $t \in \mathbb{N}$. For each density $f^{k}$, we consider two possibilities, $f_{0}^{k}$ and $f_{1}^{k}$, so that the corresponding Kullback-Leibler information numbers,  
\begin{align*}
I_1^{k} &:= \int \log \Bigl( \frac{f_1^{k}(x)}{f_0^{k}(x)} \Bigr) f_{1}^{k}(x) \, \nu^{k}(dx) \quad \text{and} \\
I_0^k &:= \int \log \Bigl(\frac{f_0^{k}(x)}{f_1^k(x)} \Bigr) f_{0}^{k}(x) \, \nu^{k}(dx),
\end{align*}
are positive and finite. The goal is to distinguish at the fusion center between the following two hypotheses:
\begin{align*} 
& \Hyp_{0}: f^{k}= f_{0}^{k}, \quad 1 \leq k \leq K \\
& \Hyp_{1}: f^{k}= f_{0}^{k}, \, k \notin \cA \quad \text{and} \quad  f^{k}= f_{1}^{k}, \, k \in \cA,
\end{align*}
where $\cA \subset \{1, \ldots, K\}$ is a subset of sensors that belongs to some class $\cP$. The interpretation is that signal is present (resp. absent) at sensor $k$ when its observations are distributed according to $f_{1}^{k}$ (resp. $f_{0}^{k}$). Thus, the null hypothesis, $\Hyp_{0}$, represents the situation in which all sensors observe noise, whereas the alternative hypothesis, $\Hyp_{1}$, corresponds to the case that signal is present in some subset of sensors.  In what follows, we denote by $\Pro_{1}^{\cA}$ and $\Exp_{1}^{\cA}$ the probability measure and the expectation, respectively, under $\Hyp_{1}$ when the 
subset of affected sensors is $\cA$, whereas the corresponding notation under $\Hyp_{0}$ will be $\Pro_{0}$ and $\Exp_{0}$.


We will be interested in the sequential version of this  hypothesis testing problem. Thus, we assume that observations at the sensors and the fusion center are acquired sequentially and we want to select the correct hypothesis at the fusion center as soon as possible. This means that the goal is to find a \textit{sequential test}, $(T,d_{T})$, which consists of an $\{\Fc_{t}\}$-stopping time $T$ and an $\Fc_{T}$-measurable random variable $d_{T}$ that takes values in $\{0,1\}$, so that  $\Hyp_{j}$ is selected on $\{d_{T}=j, T<\infty\}$,  $j=0,1$, where $\{\Fc_t\}$ is the filtration generated by the observations of \textit{all} sensors, i.e., 
$$\Fc_{t}:=\sigma(X_{s}^{k}; 1 \leq s \leq t, 1 \leq k \leq K),  \quad t \in \mathbb{N}.$$


An ideal sequential test should have small detection delay under both hypotheses, while controlling its error probabilities below prescribed levels.
Specifically,  given $\alpha, \beta \in (0,1)$ and a class of subsets of $\{1, \ldots, K\}$, $\cP$, we set 
\begin{align*}
\ccab(\cP) := \{ (T,d_{T}): & \;\Pro_{0}(d_{T}=1) \leq \alpha \\
                    \text{and} & \; \max_{\cA \in \cP} \Pro_{1}^{\cA}(d_{T}=0) \leq \beta\},
\end{align*}                             
i.e., $\ccab(\cP)$ is the class of sequential tests whose probabilities of type-I and type-II error are bounded above by $\alpha$ and $\beta$, respectively.
Then, the problem is to find a sequential test that attains
\begin{equation} \label{infima1}
\inf_{(T,d_{T}) \in \ccab(\cP)} \Exp_{0}[T] 
\end{equation}
and, \textit{for every set $\cA \in \cP$},
\begin{equation} \label{infima2}
\inf_{(T,d_{T}) \in \ccab(\cP)} \Exp_{1}^{\cA}[T].
\end{equation}
This is indeed possible when $\cP=\{\cA\}$, that is when the subset of sensors in which signal may be present is known in advance. In this case, 
$\ccab(\cP)$ reduces to 
\begin{align*}
\ccab(\cA) := \{ (T,d_{T}): & \;\Pro_{0}(d_{T}=1) \leq \alpha \\
                    \text{and} & \; \Pro_{1}^{\cA}(d_{T}=0) \leq \beta\},
\end{align*}
and from  Wald and Wolfowitz \cite{wolf} it follows that, for any $\alpha, \beta$ so that $\alpha+\beta<1$, both (\ref{infima1}) and (\ref{infima2})
are attained by Wald's  \cite{wald} Sequential Probability Ratio Test (SPRT): 
\begin{align} \label{sprt}
\begin{split}
S^{\cA} &:= \inf\{t: Z_{t}^{\cA} \notin (-A,B)\} \\
d_{S^{\cA}} &:= \left\{\begin{array}{cl}1&\text{if}~Z^{\cA}_{S^{\cA}} \geq B\\
0&\text{if}~ Z^{\cA}_{S^{\cA}} \leq -A
\end{array}\right.
\end{split}
\end{align}
where $A,B$ are positive thresholds selected so that  $\Pro_{0}( d_{S^{\cA}}=1)=\alpha$ and $\Pro_{1}^{\cA}(d_{S^{\cA}}=0)=\beta$, whereas 
$Z^{\cA}$ is the log-likelihood ratio process of $\Pro_{1}^{\cA}$ over  $\Pro_{0}$.  Since we have assumed that observations coming from different sensors are independent (an assumption that is not required for the optimality of the SPRT), it is clear that $Z_{t}^{\cA}=\sum_{k \in \cA} Z_{t}^{k}$, where 
$$Z_{t}^{k} = Z_{t-1}^{k}+ \log \frac{f_{1}^{k}(X_{t}^{k})}{f_{0}^{k}(X_{t}^{k})}; \quad Z_{0}^{k}:=0$$
is the log-likelihood ratio of the observations acquired by sensor $k$ up to  time $t$. 

While the optimality of the SPRT holds for any given $\alpha, \beta$ so that $\alpha+\beta <1$, closed-form expressions for its operating characteristics 
are, in general, available only in an asymptotic setup, that is  as $\alpha, \beta  \rightarrow 0$. In what follows, whenever $\alpha$ and $\beta$ go to 0 simultaneously, we will assume implicitly that $|\log \alpha|/|\log \beta|$ converges to some positive constant and we will write $x \sim y$ when  $\lim(x/y) = 1$ and $x \succ y$  (resp. $x \prec y$) when $\liminf(x/y) \geq 1$ (resp. $\limsup(x/y) \leq 1$). Then, it is well known that 
\begin{align} \label{sprt_perform}
\Exp_{0}[S^{\cA}] \sim \frac{|\log \beta|}{I_{0}^{\cA}}  \quad \text{and} \quad \Exp_{1}^{\cA}[S^{\cA}] &\sim \frac{|\log \alpha|}{I_{1}^{\cA}},
\end{align}
where $I_{1}^{\cA}:=\Exp_{1}^{\cA}[Z_{1}^{\cA}]$ and $I_{0}^{\cA}:=\Exp_{0}[-Z_{1}^{\cA}]$ are the  Kullback-Leibler information numbers between $\Pro_{0}$ and $\Pro_{1}^{\cA}$, which --due to the assumption of independence across sensors-- take the form $I_{j}^{\cA}=\sum_{k \in \cA} I_{j}^{k}$, $j=0,1$.

When the class $\cP$ is not a singleton, i.e., when the alternative hypothesis is composite, it is not possible to find a sequential test that 
attains (\ref{infima2}) \textit{for every subset $\cA \in \cP$}. For this reason, we need to restrict ourselves to sequential tests that are optimal in an asymptotic sense. Therefore, given a class $\cP$ of subsets of $\{1, \ldots, K\}$, we will say that a sequential test $(\tilde{T},\tilde{d}) \in \ccab(\cP)$ is \textit{asymptotically optimal} under $\Hyp_{0}$, if 
\begin{align*}
& \Exp_{0}[\tilde{T}] \sim \inf_{(T,d_{T}) \in \ccab(\cP)} \Exp_{0}[T] 
\end{align*}
and under $\Hyp_{1}$, if \textit{for  every $\cA \in \cP$}
\begin{align*}
& \Exp_{1}^{\cA}[\tilde{T}] \sim \inf_{(T,d_{T}) \in \ccab(\cP)} \Exp_{1}^{\cA}[T].
\end{align*}
A number of asymptotically optimal (under both hypotheses) sequential tests have been proposed and studied 
in the case that  signal may be present in \textit{at most one sensor}, that is when 
\begin{equation*} 
\cP=\{\cA: |\cA|=1\}=\{\{k\}, 1 \leq k \leq K\},
\end{equation*}
where $|\cA|$ represents the cardinality of $\cA$. An example of such a test is given by the SPRT-bank, according to which each sensor runs an SPRT locally,
transmits its decision to the fusion center and the latter stops and selects $\Hyp_{1}$ the first time that any sensor makes a selection in favor of the alternative, whereas it stops and selects $\Hyp_{0}$ when all sensors have made a decision in favor of the null (see, e.g., \cite{Tartakovskyetal-IEEEIT03}). 
Another asymptotically optimal sequential test 
in this setup can be obtained if $Z^{\cA}$ in (\ref{sprt}) is replaced by the generalized log-likelihood ratio statistic, $\max_{1 \leq k \leq K} Z^{k}$, or more generally  by any statistic of the form 
$$
\log \Bigl(\sum_{k=1}^{K} p_{k} e^{Z^{k}}\Bigr)  \quad \text{or} \quad 
\log \Bigl(\max_{1 \leq k \leq K} p_{k} e^{Z^{k}} \Bigr),
$$ 
where each $p_{k}$ is a positive constant \cite{Tartakovskyetal-IEEEIT03}, \cite{feltar}.

The latter approach can in principle be applied to the case that signal may be present in more than one sensors. Indeed, given \textit{any} class $\cP$, it can be shown that replacing 
$Z^{\cA}$ in (\ref{sprt}) with either 
$$
\log \Bigl(\sum_{\cB \in \cP} p_{\cB} \, e^{Z^{\cB}}\Bigr)  \quad \text{or} \quad 
\log \Bigl(\max_{\cB \in \cP} p_{\cB} \, e^{Z^{\cB}} \Bigr),
$$ 
where $Z^{\cB}:=\sum_{k \in \cB} Z^{k}$ and each $p_{\cB}$ is a positive constant, leads to an asymptotically optimal sequential test.
However, this test may not be implementable in practice, even for a moderate number of sensors. Consider, for example, the completely \textit{unstructured} case, where there is absolutely no prior information regarding the set of affected sensors and $\cP$ is given by 
\begin{equation} \label{powerset}
\cP=\{\cA: 1 \leq |\cA| \leq K\}.
\end{equation}
Then, the implementation of the above sequential tests demands summing/maximizing $2^{K}$ statistics at every time $t$, a requirement that may be prohibitive in practice. 

\subsection{Main contributions} 
In the present paper, we focus on the case that $\cP$ is given by  (\ref{powerset}), i.e., we assume that signal may be present 
in \textit{any} subset of sensors under the alternative hypothesis.
In this context, we propose a  class of sequential tests, whose implementation at any time $t$ requires $K$ (instead of $2^{K}$) operations, 
and we establish their asymptotically optimality. Specifically, we set 
$$T^{*}:=\min \{\hat{T}_{B}, \check{T}_{A}\}, \quad d^{*}:= \left\{\begin{array}{cl}1&\text{if}~\hat{T}_{B} \leq \check{T}_{A}\\
0&\text{if}~ \hat{T}_{B}> \check{T}_{A}
\end{array}\right.$$
where $\hat{T}_{B}$ and $\check{T}_{A}$ are one-sided stopping times of the form 
\begin{align*}
\check{T}_{A} &:= \min\{t: \check{Z}_{t} \leq -A \}; \quad \check{Z}_{t}:=\max_{1 \leq k  \leq K} Z^{k}_{t} \\
\hat{T}_{B} &:= \min\{t: \hat{Z}_{t} \geq B \}; \quad \hat{Z}_{t}:=\sum_{k=1}^{K} \hat{Z}^{k}_{t} ,
\end{align*}
and each $\hat{Z}^{k}$ is an $\{\Fc_{t}^{k}\}$-adapted statistic  that should be chosen appropriately. Our main contribution in this work is that we show how to select these statistics, as well as the thresholds $A$ and $B$, in order to guarantee the asymptotic optimality  of the proposed sequential test.
Thus, in Section \ref{sec2} we show that $(T^{*}, d^{*})$ is asymptotically optimal under both hypotheses, when for every $1 \leq k \leq K$ and $t \in \mathbb{N}$
\begin{equation} \label{condi}
\hat{Z}_{t}^{k} \leq M_{t}^{k}:= \max_{1 \leq s \leq t} Z_{s}^{k},
\end{equation}
there is a constant $\Di \geq 0$ so that 
\begin{equation}  \label{condi2}
\hat{Z}^{k}_{t} \geq \max\{Z^{k}_{t}-\Di, 0\} 
\end{equation}
and thresholds  $A$ and $B$ are selected so that 
\begin{equation} \label{AB}
A=A_{\beta}:= |\log \beta| \quad \text{and} \quad B=B_{\alpha}:=F^{-1}(\alpha),
\end{equation}
where $F^{-1}$ is the inverse of the survival function of the Erlang distribution  with parameters $1$ and $K$, i.e., 
\begin{equation} \label{F}
F(x):=e^{-x}  \sum_{j=0}^{K-1} \frac{x^{j}}{j!}, \quad x>0.
\end{equation} 


Conditions (\ref{condi})-(\ref{condi2}) are clearly satisfied when each $\hat{Z}^{k}_{t}$ is chosen as the positive part of $Z_{t}^{k}$, $\max\{Z_{t}^{k},0\}$, in which case $\Di=0$. In Section \ref{sec3}, we show that if each sensor $k$ communicates with the fusion center only when $Z^{k}$ increases by $\Di>0$ since the previous communication time, then selecting $\hat{Z}^{k}$ as the value of $Z^{k}$ at the most recent communication time also satisfies conditions (\ref{condi})-(\ref{condi2}). Furthermore, we show that the asymptotic optimality  of $(T^{*}, d^{*})$ remains valid  in this context, even with an asymptotically low rate of communication.

This infrequent communication is a very important property in applications characterized  by limited  bandwidth, 
where it is necessary to design schemes that require minimal transmission activity  from the sensors to the fusion center (see, e.g., \cite{tsi}, \cite{vij}).  Such communication constraints have motivated the problem of \textit{decentralized} sequential hypothesis testing (see, e.g., \cite{sama}- \cite{yil}), where each sensor is required to transmit a small number of bits whenever it communicates with the fusion center. However, 
in this literature, it is typically assumed that the set of affected sensors is known in advance (i.e., $\cP=\{\cA\}$) and asymptotically optimal decentralized sequential tests have been proposed only under this assumption (see \cite{mei}, \cite{felmoust}, \cite{yil}). Our second main contribution in the present work is that we construct a decentralized sequential  test  which requires  \textit{infrequent} transmission of \textit{one-bit messages} from the sensors and we establish its asymptotic optimality when $\cP$ is given by (\ref{powerset}).

The remaining paper is organized as follows: in Section \ref{sec2} we state and prove the main results of the paper and in Section \ref{sec3} we consider the decentralized setup. In Section \ref{sec4} we discuss certain extensions of this work, which will be presented elsewhere.

\section{Main results} \label{sec2}
In what follows, $\cP$ is given by (\ref{powerset}). We start by obtaining an asymptotic lower bound for the optimal performance under each hypothesis. 

\vspace{0.3cm}

\begin{theorem} \label{theo0}
As $\alpha, \beta \rightarrow 0$
\begin{align} \label{deli2}
\inf_{(T,d_{T}) \in \ccab(\cP)} \Exp_{0}[T] &\succ  \frac{|\log \beta|}{ \min_{1 \leq k  \leq K} I_{0}^{k}}
\end{align}
and,  for every $\cA \in \cP$,
\begin{equation} \label{deli4}
\inf_{(T,d_{T}) \in \ccab(\cP)} \Exp_{1}^{\cA}[T] \succ  \frac{|\log \alpha|}{I_{1}^{\cA}}.
\end{equation}
\end{theorem}

\vspace{0.3cm}

\begin{proof}
Since $\ccab(\cP) \subset \ccab(\cA)$ for any $\cA \in \cP$, 
\begin{align*}
\inf_{(T,d_{T}) \in \ccab(\cP)} \Exp_{1}^{\cA}[T] &\geq \inf_{(T,d_{T}) \in \ccab(\cA)} \Exp_{1}^{\cA}[T] \sim \frac{|\log \alpha|}{I_{1}^{\cA}}
\end{align*}
where the asymptotic equality follows from (\ref{sprt_perform}). This  proves (\ref{deli4}). In a similar way we can show that
\begin{align*}
\inf_{(T,d_{T}) \in \ccab(\cP)} \Exp_{0}[T] & \geq \inf_{(T,d_{T}) \in \ccab(\cA)} \Exp_{0}[T] \sim \frac{|\log \beta|}{I_{0}^{\cA}} 
\end{align*}
and  optimizing the asymptotic lower bound over $\cA \in \cP$ we obtain
\begin{align*}
\inf_{(T,d_{T}) \in \ccab(\cP)} \Exp_{0}[T] &\succ \max_{ \cA \in \cP} \frac{|\log \beta|}{I_{0}^{\cA}} = \frac{|\log \beta|}{ \min_{\cA \in \cP} I_{0}^{\cA}}.
\end{align*}
Since $I_{0}^{\cA}=\sum_{k \in \cA} I_{0}^{k}$ and $I_{0}^{k}>0$ for every $k$, it is clear that 
$\min_{\cA \in \cP} I_{0}^{\cA}= \min_{1 \leq k \leq K} I_{0}^{k}$, which proves (\ref{deli2}).
\end{proof}

\vspace{0.3cm}

In the following theorem we show that selecting $A$ and $B$ according to (\ref{AB}) guarantees that $(T^{*},d^{*}) \in \ccab(\cP)$, 
as long as each statistic $\hat{Z}^{k}$ satisfies (\ref{condi}). 

\vspace{0.3cm}

\begin{theorem}  \label{theo1}
If  $A$ and $B$ are selected according to (\ref{AB}) and each $\hat{Z}^{k}$ satisfies (\ref{condi}), then $(T^{*}, d^{*}) \in \ccab(\cP)$.
\end{theorem}

\vspace{0.2cm}

\begin{proof}
For any $A,B>0$ we have
$$\Pro_{0} (\hat{T}_{B} \leq \check{T}_{A}) \leq \Pro_{0} (\hat{T}_{B}< \infty) =\lim_{t \rightarrow \infty} \Pro_{0} (\hat{T}_{B} \leq t) $$
and for any $t\in \mathbb{N}$
\begin{align*}
\Pro_{0} (\hat{T}_{B} \leq t) &= \Pro_{0} \Bigl(\max_{0 \leq s \leq t} \sum_{k=1}^{K} \hat{Z}^{k}_{s} \geq B \Bigr) 
\leq  \Pro_{0}\Bigl(\sum_{k=1}^{K} M^{k}_{t} \geq B \Bigr),
\end{align*}
where the inequality is due to \eqref{condi}. Now, for any given $k$ and $t$, it is clear that 
\begin{align*} 
\Pro_{0}(M^{k}_{t}\geq B) &= \Pro_{0}(S^{k}_{B} \leq t) \leq \Pro_{0}(S^{k}_{B} < \infty),
\end{align*}
where $S_{B}^{k}:= \inf \{ t: Z_{t}^{k} \geq B\}$, and from Wald's likelihood ratio identity it follows that
\begin{align*}
\Pro_{0}(S^{k}_{B} < \infty) =\Exp_{1}^{k} \Bigl[e^{-Z^{k}_{S^{k}_{B}}}\Bigr] \leq e^{-B},
\end{align*}
where $\Exp^{k}_{1}$ is expectation with respect to $\Pro_{1}^{k}$, the probability  measure under which $f^{k}=f_{1}^{k}$ and $f^{j}=f_{0}^{j}$ for $j \neq k$.
The last two relationships imply that, for any given $k$ and $t$, $\Pro_{0}(M^{k}_{t}\geq B) \leq e^{-B}$, which means that the random variable  $M^{k}_{t}$ is stochastically dominated by an exponential random variable with rate 1. Since, due to the assumed independence across sensors, $M_{t}^{1},  \ldots, M_{t}^{K}$ are independent, this implies that $\sum_{k=1}^{K}M^{k}_{t}$ is stochastically dominated by an Erlang random variable with parameters 1 and $K$, i.e., 
$$\Pro_{0}\Bigl(\sum_{k=1}^{K} M^{k}_{t} \geq B \Bigr) \leq F(B),$$ where $F(x)$ is defined in (\ref{F}). From the latter observation and the definition of $B_{\alpha}$ it follows that for any $A>0$:
$$
\Pro_{0}(\hat{T}_{B_{\alpha}} \leq \check{T}_{A}) \leq F(B_{\alpha})= \alpha.
$$
Furthermore, for any given $\cA  \in \cP$, from Wald's likelihood ratio identity it follows that for any $A,B>0$
\begin{align*} 
\Pro_{1}^{\cA} (\check{T}_{A} < \hat{T}_{B}) &\leq  \Pro_{1}^{\cA} (\check{T}_{A} < \infty) \\
&= \Exp_{0}\Bigl[ e^{\sum_{k \in \cA} Z^{k}_{\check{T}_{A}}}] \\
&\leq e^{-|\cA| A} \leq e^{-A},   
\end{align*}
where the second inequality holds because $Z^{k}_{\check{T}_{A}}\leq -A$ on $\{\check{T}_{A} < \infty\}$ for every $1 \leq k \leq K$ and the third one 
because $|\cA|  \geq 1$ for any $\cA \in \cP$.   Consequently,  
\begin{align*} 
\max_{\cA \in \cP} \Pro_{1}^{\cA} (\check{T}_{A} < \hat{T}_{B}) &\leq  e^{-A}
\end{align*}
and from  the definition of $A_{\beta}$ it follows that for any $B>0$
\begin{align*} 
\max_{\cA \in \cP} \Pro_{1}^{\cA} (\check{T}_{A_{\beta}} < \hat{T}_{B}) &\leq  e^{-A_{\beta}}=\beta,  
\end{align*}
which completes the proof.
\end{proof}

\vspace{0.3cm}

In the following theorem we show that if $A$ and $B$ are selected according to (\ref{AB}) and each  statistic $\hat{Z}^{k}$ satisfies (\ref{condi2}), then
$(T^{*},d^{*})$ attains the asymptotic lower bounds in (\ref{deli2}) and (\ref{deli4}). 

\vspace{0.3cm}

\begin{theorem} \label{theo2}
(i)  As $A \rightarrow \infty$,
\begin{equation} \label{deli1}
\Exp_{0}[T^{*}] \prec \frac{A}{\min_{1 \leq k \leq K} I_{0}^{k}}
\end{equation}
and $(T^{*},d^{*})$ attains the asymptotic lower bound in (\ref{deli2}) when $A=A_{\beta}$.

(ii) If each $\hat{Z}^{k}$ satisfies (\ref{condi2}), then as $B \rightarrow \infty$
\begin{equation} \label{deli3}
\Exp_{1}^{\cA}[T^{*}] \prec \frac{B+ \sum_{k \in \cA} \Di}{I_{1}^{\cA}} \quad \text{for every} \; \cA \in \cP
\end{equation} 
and  $(T^{*},d^{*})$ attains the asymptotic lower bound in (\ref{deli1}) when $B=B_{\alpha}$.
\end{theorem}

\vspace{0.3cm}

\begin{proof}
The proof of (i) is a direct consequence of Theorem 2 in \cite{Tartakovskyetal-IEEEIT03}.
In order to prove (ii), we observe that for any $k$ and $t$
\begin{align*}
\sum_{k=1}^{K} \hat{Z}^{k}_{t} &= \sum_{k \in \cA} \hat{Z}^{k}_{t}  + \sum_{k \notin \cA} \hat{Z}^{k}_{t}  \\
                               &\geq  \sum_{k \in \cA} (Z^{k}_{t}-\Di)  = Z_{t}^{\cA}-  \sum_{k \in \cA} \Di  ,
\end{align*}
where the inequality is due to (\ref{condi2}). As a result, $$T^{*} \leq \hat{T}_{B} \leq \inf\Bigl\{t: Z_{t}^{\cA}\geq B+ \sum_{k \in \cA} \Di\Bigr\}$$
and taking expectations we obtain (\ref{deli3}). From this relationship and (\ref{deli4}) it is clear that it 
suffices to show that $B_{\alpha} \sim |\log \alpha|$ as $\alpha \rightarrow 0$. Indeed, 
taking logarithms in the definition of $B_{\alpha}$ in (\ref{AB})-(\ref{F}) we have 
\begin{align} \label{ba}
|\log \alpha|  &= B_{\alpha} - \log  \Bigl(\sum_{j=0}^{K-1} \frac{B_{\alpha}^{j}}{j!} \Bigr) \sim B_{\alpha},
\end{align} 
which completes the proof. 
\end{proof}

\vspace{0.3cm}

From the previous theorems it follows that selecting $A$ and $B$ according to (\ref{AB}) and the statistics $\{\hat{Z}^{k}, 1\leq k \leq K\}$ so that  (\ref{condi})-(\ref{condi2}) hold guarantees the asymptotic optimality of $(T^{*},d^{*})$ under both hypotheses, when $\cP$ is given by (\ref{powerset}). 
Let us add a few remarks to this statement:

\begin{enumerate}
\item Conditions (\ref{condi})-(\ref{condi2}) are clearly satisfied when $\hat{Z}^{k}_{t}= \max\{Z_{t}^{k},0\}$. An alternative specification that satisfies these conditions is presented in the next section. 

\item Condition (\ref{condi2}) is \textit{not} needed for $(T^{*},d^{*})$ to belong in $\ccab(\cP)$ and to be asymptotically optimal under $\Hyp_{0}$.

\item The asymptotic optimality of $(T^{*},d^{*})$ remains valid even if $\Di\rightarrow \infty$ for one or more $k$, as long as $\Di=o(|\log \alpha|)$. In the next section we will show that, with a particular specification for $\hat{Z}^{k}$, this property has an interesting interpretation in terms of the communication requirements of the proposed scheme.

\end{enumerate}

\section{The decentralized setup} \label{sec3}
Let us first note that the one-sided sequential test $\check{T}_{A}$ is an one-shot scheme; it requires that each sensor communicate with the fusion center
at most once, as soon as its local log-likelihood statistic takes a value smaller than $-A$, at which time it simply needs to transmit a one-bit message to the fusion center, informing it about this development. 

On the other hand, the implementation of the stopping rule $\hat{T}_{B}$ can be much more demanding from a communication point of view. For example, if we set $\hat{Z}^{k}_{t}= \max\{Z_{t}^{k},0\}$, sensor $k$ needs to transmit the actual value of $Z^{k}$ at every time $t$ (or at least whenever it is positive). As we discussed in the Introduction, this may not be possible in applications characterized by  bandwidth constraints. 

In what follows, we assume that thresholds $A$ and $B$ are selected according to (\ref{AB}) and  our goal is to suggest specifications for $\{\hat{Z}^{k}\}_{1 \leq k \leq K}$ that induce low transmission activity, while preserving the asymptotic optimality of the sequential test $(T^{*},d^{*})$. 

In order to achieve this, we require that each sensor $k$ communicate with the fusion center 
only at an increasing sequence of $\{\cFt^{k}\}$-stopping times, $(\tau_{n}^{k})_{n \in \mathbb{N}}$,  which are finite under $\Pro_{1}^{k}$. In other words, 
each sensor should communicate with the fusion center only at some particular time instances and, at any given time, the decision to communicate or not should 
depend exclusively on the observations that have been acquired locally at the sensor until this time.

Given such a sequence of communication times, we denote by $\tau^{k}(t)$ the instance of the most recent transmission and by $N_{t}^{k}$ the number of transmitted messages up to time $t$,  i.e., 
\begin{align*}
\tau^{k}(t) &:= \max\{\tau_{n}^{k}:\tau_{n}^{k} \leq t\}, \quad N_{t}^{k} := \max\{n: \tau_{n}^{k} \leq t\}.
\end{align*}
At any given time $t$, the value of $Z^{k}$ at the most recent communication instance, 
\begin{equation} \label{second}
Z^{k}_{\tau^{k}(t)}= \sum_{n=1}^{N_{t}^{k}} \ell_{n}^{k}, \quad  \ell_{n}^{k}:= Z_{\tau_{n}^{k}}^{k}- Z_{\tau_{n-1}^{k}}^{k},
\end{equation}
cannot be larger than $M_{t}^{k}$, the maximum value of $Z^{k}$  up to time $t$. Indeed, note that $Z^{k}_{\tau^{k}(t)}$  coincides with $Z^{k}$ at the communication times $(\tau_{n}^{k})_{n \in \mathbb{N}}$  and stays flat in between. Therefore, selecting $\hat{Z}_{t}^{k}$  according to (\ref{second}) 
satisfies condition (\ref{condi}) and, consequently, it guarantees that $(T^{*},d^{*})$ belongs in $\ccab(\cP)$ and is asymptotically optimal under $\Hyp_{0}$.


When, in particular,  the communication times are described by the recursion
\begin{align} \label{tau}
\tau_{n}^{k} &:=\inf \{t \geq \tau_{n-1}^{k}: Z_{t}^{k}- Z_{\tau_{n-1}^{k}}^{k} \geq \Di\}, \quad n \in \mathbb{N} 
\end{align}
where $\tau_{0}^{k}:=0$ and  $\Di$ is a positive constant, then it is straightforward to see that, for any time $t$, $Z^{k}_{\tau^{k}(t)}\geq 0$
and $Z_{t}^{k}- Z^{k}_{\tau^{k}(t)} \leq \Di$.  Therefore, in the case of the communication scheme (\ref{tau}),  setting $\hat{Z}^{k}_{t}$ equal to $Z^{k}_{\tau^{k}(t)}$ satisfies condition (\ref{condi2}) as well and implies that  $(T^{*},d^{*})$ is  asymptotically optimal also under $\Hyp_{1}$. Furthermore, the final remark in the end of Section \ref{sec2} suggests that the latter asymptotic optimality property is preserved even with an asymptotically low rate of communication from one or more sensors, as the constant $\Di$ in this setup controls the average period of communication at sensor $k$.

From the right-hand side in (\ref{second}) it is clear that selecting $\hat{Z}^{k}_{t}$ as $Z^{k}_{\tau^{k}(t)}$ requires
that at each time  $\tau_{n}^{k}$ sensor $k$ transmit to the   fusion center (with an ``infinite-bit'' message) the exact value of $\ell_{n}^{k}$, the ``realized'' local log-likelihood ratio between $\tau_{n-1}^{k}$ and $\tau_{n}^{k}$. However, if one insists that a small number of bits be transmitted at each communication, which is the main requirement in decentralized sequential testing \cite{veer}, then this selection is no longer acceptable. Nevertheless, in the case of the communication scheme (\ref{tau}), it is intuitively clear that the value of each $\ell_{n}^{k}$ should be close to $\Di$, at least when $Z^{k}$ does not have ``heavy tails'' and/or  $\Di$ is ``large''. This suggests selecting each $\hat{Z}^{k}_{t}$ according 
\begin{align} \label{third}
 \sum_{n=1}^{N_{t}^{k}} \Di = \Di N_{t}^{k},
\end{align}
a selection that requires transmission of a \textit{single bit} from each sensor at each communication time. Moreover, for every time $t$ it is clear that 
\begin{equation} \label{ineqa2}
\D \, N_{t}^{k}  \leq Z^{k}_{\tau^{k}(t)} \leq M_{t}^{k},
\end{equation} 
therefore, selecting $\hat{Z}^{k}_{t}$ according to (\ref{third}) satisfies condition (\ref{condi}) and, consequently, it
guarantees that $(T^{*},d^{*})$ belongs in $\ccab(\cP)$ and is asymptotically optimal under $\Hyp_{0}$.  
On the other hand, for every $t$ we have
\begin{equation} \label{ineqa}
Z_{t}^{k}- \Di N_{t}^{k} = Z_{t}^{k}- Z^{k}_{\tau^{k}(t)} + \sum_{n=1}^{N_{t}^{k}} \eta_{n}^{k} \leq \Di + \sum_{n=1}^{N_{t}^{k}} \eta_{n}^{k},
\end{equation}
where $\eta_{n}^{k}:= Z_{\tau_{n}^{k}}^{k}- Z_{\tau_{n-1}^{k}}^{k}-\Di$ is the random, non-negative overshoot associated with the $n^{th}$ transmission from sensor $k$. This means that selecting each $\hat{Z}^{k}$ according to (\ref{third}) does not satisfy condition (\ref{condi2}), therefore Theorem \ref{theo2}(ii)
can no longer be applied to establish the asymptotic optimality of $(T^{*},d^{*})$ under $\Hyp_{1}$. 

Nevertheless, in the following theorem  we show that this property remains valid if two additional conditions are satisfied. 
The first is that, for every $k \in \cA$, each $Z_{1}^{k}$ must have a finite second moment under $\Pro_{1}^{k}$, i.e., 
\begin{equation} \label{moment}
\Exp_{1}^{k}[(Z_{1}^{k})^{2}]= \int \log \Bigl( \frac{f_1^{k}(x)}{f_0^{k}(x)} \Bigr)^{2} \, f_{1}^{k}(x) \, \nu^{k}(dx) <\infty,
\end{equation} 
a condition that guarantees that
\begin{equation} \label{C}
C_{\cA}:=\max_{k \in \cA} \, \sup_{\Di>0} \, \Exp_{1}^{k}[\eta^{k}_{1}]
\end{equation}
is a finite quantity for any given $\{\Di, k \in \cA\}$ and an $\calo(1)$ term as $\Di \rightarrow \infty$ for every $k \in \cA$. 
 The second  is that, now, we \textit{must} let $\Di \rightarrow \infty$ so that $\Di=o(|\log \alpha|)$ for every $1 \leq k \leq K$,
so that each sensor does not communicate with the fusion center very frequently 
and the (unobserved) overshoots do not accumulate very fast.  

In what follows, we  denote by $\calo(\oD)$ a term that is bounded above when divided by $\oD$ as $\uD \rightarrow \infty$, where
$$\uD:=\min_{1 \leq k \leq K} \Di, \quad \oD:=\max_{1 \leq k \leq K} \Di.$$

\vspace{0.3cm}

\begin{theorem} \label{theo4}
Suppose that each sensor $k$ communicates with the fusion center at the sequence of times described by (\ref{tau}) and 
that each $\hat{Z}^{k}_{t}$ is selected according to (\ref{third}).

(i) If $A=A_{\beta}$  and   $B=B_{\alpha}$, then $(T^{*},d^{*})$ belongs to $\ccab(\cP)$ and 
attains the asymptotic lower bound in (\ref{deli2}). 

(ii) If (\ref{moment}) holds, then for any $B$ and $\{\Di, 1 \leq k \leq K\}$
\begin{align} \label{deli5}
\Exp_{1}^{\cA}[T^{*}] &\leq \frac{1}{I_{1}^{\cA}} \Bigl[ \calo(\oD) + \Bigl(1+\frac{C_{\cA}}{\uD}\Bigr) B \Bigr],
\end{align}
and $(T^{*},d^{*})$ attains the asymptotic lower bound in (\ref{deli4}) when $B=B_{\alpha}$, as long as $\Di \rightarrow \infty$ so that $\Di=o(|\log \alpha|)$ for every $1 \leq k \leq K$.
\end{theorem}

\vspace{0.3cm}

\begin{proof} 
The proof of (i) follows from (\ref{ineqa2}) and Theorems  \ref{theo0}(i), \ref{theo1} and \ref{theo2}(i). In order to prove (ii), 
we start with the observation that $T^{*} \leq \hat{T}_{B}$  for any thresholds $A,B$ and that for  any subset $\cA$ we have
\begin{align} \label{repli}
\begin{split}
I_{1}^{\cA} \, \Exp_{1}^{\cA}[\hat{T}_{B}] &= \Exp_{1}^{\cA}[Z_{\hat{T}_{B}}^{\cA}]  \\
&=  \sum_{k \in \cA} \Exp_{1}^{\cA}[(Z^{k}-\hat{Z}^{k})_{\hat{T}_{B}}] +   \sum_{k \in \cA} \Exp_{1}^{\cA}[\hat{Z}^{k}_{\hat{T}_{B}}] \\
&\leq \sum_{k \in \cA} \Exp_{1}^{\cA}[(Z^{k}-\hat{Z}^{k})_{\hat{T}_{B}}] + \Exp_{1}^{\cA}[\hat{Z}_{\hat{T}_{B}}],
\end{split}
\end{align}
where the equality follows from Wald's identity, whereas the inequality holds because $\sum_{k \in \cA} \hat{Z}^{k} \leq \hat{Z}$  
whenever every $\hat{Z}^{k}$ is non-negative, as it is the case with (\ref{third}).

For any $k \in \cA$, setting $t=\hat{T}_{B}$ in (\ref{ineqa}) and strengthening the inequality we have 
\begin{align} \label{des}
(Z^{k}-\hat{Z}^{k})_{\hat{T}_{B}} &\leq \Di + \sum_{n=1}^{N_{\hat{T}_{B}}^{k}+1} \eta_{n}^{k}.
\end{align}
Moreover, setting $\Gc_{n}^{k}:=\Fc_{\tau_{n}^{k}}$, $n \in \mathbb{N}$, we can see that $N_{\hat{T}_{B}}^{k}+1$ is a $\Pro_{1}^{\cA}$-integrable,
$\{\Gc_{n}^{k}\}_{n \in \mathbb{N}}$-adapted stopping time and $(\eta_{n}^{k})_{n \in \mathbb{N}}$  a sequence of $\{\Gc_{n}^{k}\}$-adapted, i.i.d. random variables with finite expectation, $\Exp_{1}^{\cA}[\eta_{1}^{k}]=\Exp_{1}^{k}[\eta_{1}^{k}]$.  As a result, from Wald's first identity it follows that 
for every $k \in \cA$:
$$\Exp_{1}^{\cA}\Bigl[\sum_{n=1}^{N_{\hat{T}_{B}}^{k}+1} \eta_{n}^{k}\Bigr]= \Exp_{1}^{\cA}[N_{\hat{T}_{B}}^{k}+1] \, \Exp_{1}^{k}[\eta_{1}^{k}].$$
Therefore,  taking expectations in (\ref{des}) and recalling the definition of $C_{\cA}$ in (\ref{C}) we have
\begin{align*}
\Exp_{1}^{\cA}[(Z^{k}-\hat{Z}^{k})_{\hat{T}_{B}}]  &\leq \Di + C_{\cA} +  C_{\cA}  \, \Exp_{1}^{\cA}[N_{\hat{T}_{B}}^{k}].
\end{align*}
Then, summing over $k \in \cA$  and setting $N_{t}:=\sum_{k=1}^{K} N_{t}^{K}$, we obtain
\begin{align*}
\sum_{k \in \cA} \Exp_{1}^{\cA}[(Z^{k}-\hat{Z}^{k})_{\hat{T}_{B}}]
&\leq  \sum_{k \in \cA} [\Di+C_{\cA}]  + C_{\cA}  \, \Exp_{1}^{\cA}[ N_{\hat{T}_{B}}].
\end{align*}
However, since each $\hat{Z}^{k}$ is selected according to (\ref{third}), then it is clear  that  
$\hat{Z}_{t}  \geq \uD \, N_{t}$ for  every $t$. Therefore,  
\begin{align*}
\sum_{k \in \cA} \Exp_{1}^{\cA}[(Z^{k}-\hat{Z}^{k})_{\hat{T}_{B}}]
&\leq  \sum_{k \in \cA} [\Di+C_{\cA}]  + C_{\cA} \, \frac{ \Exp_{1}^{\cA}[\hat{Z}_{\hat{T}_{B}}]}{\uD}
\end{align*}
and from  (\ref{repli}) it follows that $I_{1}^{\cA} \, \Exp_{1}^{\cA}[T^{*}]$ is bounded above by
\begin{align} \label{uppers}
\sum_{k \in \cA} [\Di+C_{\cA}] + \Bigl(1+ \frac{C_{\cA}}{\uD}\Bigr) \, \Exp_{1}^{\cA}[\hat{Z}_{\hat{T}_{B}}].
\end{align} 
But since each $\hat{Z}^{k}$ is selected according to (\ref{third}), it is clear that the overshoot $\hat{Z}_{\hat{T}_{B}} - B$ cannot take a value larger than 
$\sum_{k=1}^{K}\Di \leq K \oD$. As a result, $\hat{Z}_{\hat{T}_{B}}\leq B + K \oD$ and the upper bound (\ref{uppers}) takes the form 
$$
\Bigl[\sum_{k \in \cA} (\Di+C_{\cA}) + \Bigl(1+ \frac{C_{\cA}}{\uD}\Bigr) K \oD \Bigr] + \Bigl(1+ \frac{C_{\cA}}{\uD}\Bigr) B.
$$
Then, in order to prove (\ref{deli5}), it suffices to note that the first two terms in the latter expression are $\calo(\oD)$,
since  $C_{\cA}$ is an $\calo(1)$ term as $\uD \rightarrow \infty$, due to assumption (\ref{moment}).

Finally, since $B_{\alpha} \sim |\log \alpha|$ as $\alpha \rightarrow 0$ (recall \eqref{ba}), 
from (\ref{deli5}) it follows that $(T^{*},d^{*})$ attains the asymptotic lower bound in (\ref{deli4}) when $B=B_{\alpha}$, 
as long as $\uD \rightarrow \infty$ so that $\oD=o(|\log \alpha|)$, which completes the proof. 
\end{proof}

\section{Extensions} \label{sec4}
Theorems \ref{theo0},  \ref{theo1} and \ref{theo2} do not rely heavily on the assumed i.i.d. structure of the sensor observations. Thus, it can be shown that the asymptotic optimality of $(T^{*},d^{*})$ remains valid for any statistical model (in discrete or continuous time) that preserves the asymptotic optimality of the SPRT. Moreover, the above results can be generalized in the case that a lower bound,  $\uK \geq 1$, and an upper bound, $\oK \leq K$, are available on the number of affected sensors, that is when  $\cP=\{\cA: \uK \leq |\cA| \leq \oK\}$.

Furthermore, it is straightforward to generalize the decentralized sequential test described in Section \ref{sec3}, so that more than one bits are transmitted per communication. These additional bits can be utilized for the quantization of the unobserved overshoots and can improve the performance of the proposed test in the case of high rates of communication. 

Finally, we should note that all these extensions, which will be presented elsewhere, require the assumption of independence across sensors. Removing this assumption remains an open problem.

\section*{Acknowledgments}
This work was supported by the U.S.\ Air Force Office of Scientific Research under MURI grant FA9550-10-1-0569, by the U.S.\ Defense Threat Reduction Agency under grant HDTRA1-10-1-0086, by the U.S.\  Defense Advanced Research Projects Agency under grant W911NF-12-1-0034, by the U.S.\ National Science Foundation under grants CCF-0830419, EFRI-1025043, and DMS-1221888 and by the U.S. Army Research Office under grant W911NF-13-1-0073 at the University of Southern California, Department of Mathematics.



\bibliographystyle{IEEEtran}

\end{document}